\documentclass[journal,onecolumn,12pt]{article}

\usepackage{amssymb}
\usepackage{caption2}
\usepackage{amsmath}
\usepackage{multirow}
\usepackage{amsthm}
\usepackage{graphicx}
\usepackage{CJK}
\usepackage{enumerate}
\usepackage{appendix}
\usepackage{bm}
\usepackage{tikz}
\usetikzlibrary{positioning,chains,fit,shapes,calc}

\newtheorem{theorem}{Theorem}[section]

\newtheorem{definition}[theorem]{Definition}

\newtheorem{remark}[theorem]{Remark}

\begin{document}
\title{Some results on similar configurations in subsets of $\mathbb{F}_q^d$}

\author{Chengfei Xie\thanks{C. Xie is with the School of Mathematical Sciences, Capital Normal University, Beijing 100048, China (e-mail: cfxie@cnu.edu.cn).}
~and Gennian Ge\thanks{Corresponding author. G. Ge is with the School of Mathematical Sciences, Capital Normal University, Beijing 100048, China (e-mail: gnge@zju.edu.cn). The research of G. Ge was supported by the National Key Research and Development Program of China under Grant 2020YFA0712100 and Grant 2018YFA0704703, the National Natural Science Foundation of China under Grant 11971325  and Grant 12231014,  and Beijing Scholars Program.}}

\maketitle

\begin{abstract}
In this paper, we study problems about the similar configurations in $\mathbb{F}_q^d$. Let $G=(V, E)$ be a graph, where $V=\{1, 2, \ldots, n\}$ and $E\subseteq{V\choose2}$. For a set $\mathcal{E}$ in $\mathbb{F}_q^d$, we say that $\mathcal{E}$ contains a pair of $G$ with dilation ratio $r$ if there exist distinct $\bm{x}_1, \bm{x}_2, \ldots, \bm{x}_n\in\mathcal{E}$ and distinct $\bm{y}_1, \bm{y}_2, \ldots, \bm{y}_n\in\mathcal{E}$ such that $\|\bm{y}_i-\bm{y}_{j}\|=r\|\bm{x}_i-\bm{x}_j\|\neq0$ whenever $\{i, j\}\in E$, where $\|\bm{x}\|:=x_1^2+x_2^2+\cdots+x_d^2$ for $\bm{x}=(x_1, x_2, \ldots, x_d)\in\mathbb{F}_q^d$. We show that if $\mathcal{E}$ has size at least $C_kq^{d/2}$, then $\mathcal{E}$ contains a pair of $k$-stars with dilation ratio $r$, and that if $\mathcal{E}$ has size at least $C\cdot\min\left\{q^{(2d+1)/3}, \max\left\{q^3, q^{d/2}\right\}\right\}$, then $\mathcal{E}$ contains a pair of $4$-paths with dilation ratio $r$. Our method is based on enumerative combinatorics and graph theory.
\smallskip
\end{abstract}
\medskip


\smallskip


\section{Introduction}
A lot of problems in discrete geometry ask the minimum size of a set to ensure the existence of a certain structure. One of the problems of this type is the Falconer's conjecture, which states that if the Hausdorff dimension of a set $\mathcal{E}\subseteq\mathbb{R}^d$ is strictly larger than $d/2$, then the set of Euclidean distances between pairs of points in $\mathcal{E}$ has positive Lebesgue measure. We refer the readers to \cite{MR4201782}, \cite{MR4297185} and \cite{MR4055179} for the best known results.

Let $\mathbb{F}_q$ be a field with $q$ elements. In the setting of finite fields, the Erd\H{o}s-Falconer distance problem asks the minimum size of $\mathcal{E}\subseteq\mathbb{F}_q^d$ to ensure that the distance set
$$
\Delta(\mathcal{E}):=\{\|\bm{x}-\bm{y}\|:\bm{x}, \bm{y}\in\mathcal{E}\}
$$
contains a positive proportion of $\mathbb{F}_q$, where $\|\bm{x}\|:=x_1^2+x_2^2+\cdots+x_d^2$ for $\bm{x}=(x_1, x_2, \ldots, x_d)\in\mathbb{F}_q^d$. In \cite{MR2336319}, Iosevich and Rudnev showed that if $|\mathcal{E}|\geq Cq^{(d+1)/2}$ for some constant $C$, then $\Delta(\mathcal{E})=\mathbb{F}_q$. They also proved that $|\mathcal{E}|\geq Cq^{d/2}$ is a necessary condition in general. Hart et al. \cite{MR2775806} showed that the exponent $(d+1)/2$ is the best possible in odd dimensions. In \cite{MR2917133}, Chapman et al. proved that if $\mathcal{E}\subseteq\mathbb{F}_q^2$ with $q\equiv3$ {\rm (mod 4)} and $|\mathcal{E}|\geq q^{4/3}$, then $|\Delta(\mathcal{E})|\geq cq$, which improved the exponent $(d+1)/2$ in dimension 2. Bennett et al. \cite{MR3592595} extended this result to two-dimensional vector spaces over arbitrary finite fields. It is still open that whether the exponent $d/2$ is sufficient in even dimensions.

Iosevich et al. \cite{MR3959878} considered the quotient set
$$
\frac{\Delta(\mathcal{E})}{\Delta(\mathcal{E})}:=\left\{\frac{a}{b}:a\in\Delta(\mathcal{E}), b\in\Delta(\mathcal{E})\setminus\{0\}\right\}
$$
of the distance set, and obtained the following theorems.
\begin{theorem}[\cite{MR3959878}]\label{shangji}
  Let $\mathcal{E}\subseteq\mathbb{F}_q^d$ and $d\geq2$ be even. If $|\mathcal{E}|\geq9q^{d/2}$, then
  $$
\frac{\Delta(\mathcal{E})}{\Delta(\mathcal{E})}=\mathbb{F}_q.
$$
\end{theorem}
\begin{theorem}[\cite{MR3959878}]\label{shangji2}
  Let $\mathcal{E}\subseteq\mathbb{F}_q^d$ and $d\geq3$ be odd. If $|\mathcal{E}|\geq6q^{d/2}$, then
  $$
\frac{\Delta(\mathcal{E})}{\Delta(\mathcal{E})}\supseteq\mathbb{F}_q^+\cup\{0\},
$$
where $\mathbb{F}_q^+$ denotes the set of nonzero quadratic residues in $\mathbb{F}_q$.
\end{theorem}

The conclusion of Theorem \ref{shangji} means that for every $r\in\mathbb{F}_q^*$, there exist $\bm{x}, \bm{y}, \bm{x}', \bm{y}'\in\mathcal{E}$ such that $\|\bm{y}'-\bm{x}'\|=r\|\bm{y}-\bm{x}\|\neq0$. In other words, $K_{|\mathcal{E}|}$ contains a pair of edges with dilation ratio $r$. To be precise, we give the following definition.
\begin{definition}\label{zhudingyi}
  Let $G=(V, E)$ be a graph, where $V=\{1, 2, \ldots, n\}$ and $E\subseteq{V\choose2}$. For a set $\mathcal{E}$ in $\mathbb{F}_q^d$, we say that $\mathcal{E}$ contains a pair of $G$ with dilation ratio $r$ if there exist distinct $\bm{x}_1, \bm{x}_2, \ldots, \bm{x}_n\in\mathcal{E}$ and distinct $\bm{y}_1, \bm{y}_2, \ldots, \bm{y}_n\in\mathcal{E}$ such that $\|\bm{y}_i-\bm{y}_{j}\|=r\|\bm{x}_i-\bm{x}_j\|\neq0$ whenever $\{i, j\}\in E$.
\end{definition}
\begin{remark}
  When $d=2$ and $-1$ is not a quadratic residue of $\mathbb{F}_q$, $\bm{x}_i\neq\bm{x}_j$ implies $\|\bm{x}_i-\bm{x}_j\|\neq0$. In other cases, it is possible that $\bm{x}_i\neq\bm{x}_j$ and $\|\bm{x}_i-\bm{x}_j\|=0$. In order to avoid the trivial case, we add the condition that $\|\bm{x}_i-\bm{x}_j\|\neq0$ in Definition \ref{zhudingyi}.
\end{remark}
Recently, Rakhmonov \cite{2022arXiv220811579R} considered the minimum size condition to ensure that $\mathcal{E}$ contains {a pair of $G$ with dilation ratio $r$} for a certain graph $G$. For convenience, we give the basic definitions of some classes of graphs. For a natural number $n$, we write $[n]$ to denote the set $\{1, 2, \ldots, n\}$.
\begin{definition}
   Let $G=(V, E)$ be a graph, where $V=[n]$ and $E\subseteq{[n]\choose2}$.
   \begin{itemize}
     \item We say that $G$ is a $k$-path if $n=k+1$ and $E=\{\{1, 2\}, \{2, 3\}, \ldots, \{k, k+1\}\}$.
     \item We say that $G$ is a $k$-star if $n=k+1$ and $E=\{\{1, 2\}, \{1, 3\}, \ldots, \{1, k+1\}\}$.
     \item We say that $G$ is a $k$-cycle if $n=k$ and $E=\{\{1, 2\}, \{2, 3\}, \ldots, \{k-1, k\}, \{k, 1\}\}$.
   \end{itemize}
\end{definition}
Let $\mathbb{F}_q^*$ be the set of nonzero elements in $\mathbb{F}_q$. Rakhmonov obtained the following results.
\begin{theorem}[\cite{2022arXiv220811579R}]
  If $r\in\mathbb{F}_p^*$, $p$ is a prime such that $p\equiv3$ {\rm (mod 4)} and $\mathcal{E}\subseteq\mathbb{F}_p^2$ with $|\mathcal{E}|>(\sqrt3+1)p$, then  $\mathcal{E}$ contains a pair of $2$-paths with dilation ratio $r$.
\end{theorem}
\begin{theorem}[\cite{2022arXiv220811579R}]
  If $r\in\mathbb{F}_p^*$, $p$ is a prime such that $p\equiv3$ {\rm (mod 4)} and $\mathcal{E}\subseteq\mathbb{F}_p^2$ with $|\mathcal{E}|>4\sqrt3p^{3/2}$, then  $\mathcal{E}$ contains a pair of $4$-cycles with dilation ratio $r$.
\end{theorem}

In this paper, we focus on $k$-stars and $4$-paths  and we have the following theorems.
\begin{theorem}\label{diyidingli}
   Let $q$ be an odd prime power, $\mathcal{E}\subseteq\mathbb{F}_q^d$ and $k\geq2$ be an integer.
     \begin{itemize}
    \item If $q\geq5$, $d\geq2$ is even, $r\in\mathbb{F}_q^*$ and $\mathcal{E}$ has size at least $\left(31+10{k\choose2}\right)q^{d/2}$, then $\mathcal{E}$ contains a pair of $k$-stars with dilation ratio $r$.
    \item If $d\geq3$ is odd, $r\in\mathbb{F}_q^+$ and $\mathcal{E}$ has size at least $\left(4+\sqrt3{k\choose2}\right)q^{d/2}$, then $\mathcal{E}$ contains a pair of $k$-stars with dilation ratio $r$.
  \end{itemize}
\end{theorem}
\begin{theorem}\label{dierdingli}
Let $q$ be an odd prime power and $\mathcal{E}\subseteq\mathbb{F}_q^d$.
     \begin{itemize}
    \item If $q\geq5$, $2\leq d\leq4$ is even, $r\in\mathbb{F}_q^*$ and $\mathcal{E}$ has size at least $36q^{(2d+1)/3}$, then $\mathcal{E}$ contains a pair of $4$-paths with dilation ratio $r$.
    \item If $d=3$, $r\in\mathbb{F}_q^+$ and $\mathcal{E}$ has size at least $9q^{(2d+1)/3}$, then $\mathcal{E}$ contains a pair of $4$-paths with dilation ratio $r$.
        \item If $d=5$, $r\in\mathbb{F}_q^+$ and $\mathcal{E}$ has size at least $12q^{3}$, then $\mathcal{E}$ contains a pair of $4$-paths with dilation ratio $r$.
        \item If $q\geq5$, $d\geq6$ is even, $r\in\mathbb{F}_q^*$ and $\mathcal{E}$ has size at least $313q^{d/2}$, then $\mathcal{E}$ contains a pair of $4$-paths with dilation ratio $r$.
        \item If $d\geq7$ is odd, $r\in\mathbb{F}_q^+$ and $\mathcal{E}$ has size at least $313q^{d/2}$, then $\mathcal{E}$ contains a pair of $4$-paths with dilation ratio $r$.
  \end{itemize}
\end{theorem}
\begin{remark}
  We do not attempt to optimize the constant factors in Theorems \ref{diyidingli} and \ref{dierdingli}.
\end{remark}
Our idea is the following. Recall that in Definition \ref{zhudingyi}, $\bm{x}_1, \bm{x}_2, \ldots, \bm{x}_n$ are distinct and so are $\bm{y}_1, \bm{y}_2, \ldots, \bm{y}_n$. If we allow some $\bm{x}_i$ (or $\bm{y}_i$) to be equal, then there is a \textit{degenerate} pair of graphs. We have
$$
\text{the number of \textit{nondegenerate} pairs}=\text{the number of all pairs}-\text{the number of \textit{degenerate} pairs}.
$$
When the size of $\mathcal{E}$ is large, the number of all pairs is large and the number of \textit{degenerate} pairs is small. Thus $\mathcal{E}$ will contain a pair of $G$ with dilation ratio $r$. We also use tools from graph theory.

\section{Proof of Theorem \ref{diyidingli}}
In this section, we investigate the minimum size of  $\mathcal{E}$ to ensure that $\mathcal{E}$ contains a pair of $k$-stars with dilation ratio $r$ and prove Theorem \ref{diyidingli}.

Let
$$
S_k(r)=\{(\bm{x}_1, \bm{x}_2, \ldots, \bm{x}_{k+1}, \bm{y}_1, \bm{y}_2, \ldots, \bm{y}_{k+1})\in\mathcal{E}^{2k+2}:\|\bm{y}_i-\bm{y}_1\|=r\|\bm{x}_i-\bm{x}_1\|\neq0, i=2, 3, \ldots, k+1\}.
$$
Let
$$
A_{ij}=\{(\bm{x}_1, \bm{x}_2, \ldots, \bm{x}_{k+1}, \bm{y}_1, \bm{y}_2, \ldots, \bm{y}_{k+1})\in S_k(r):\bm{x}_i=\bm{x}_j\},
$$
$$
A'_{ij}=\{(\bm{x}_1, \bm{x}_2, \ldots, \bm{x}_{k+1}, \bm{y}_1, \bm{y}_2, \ldots, \bm{y}_{k+1})\in S_k(r):\bm{y}_i=\bm{y}_j\},
$$
and
$$
B=\{(\bm{x}_1, \bm{x}_2, \ldots, \bm{x}_{k+1}, \bm{y}_1, \bm{y}_2, \ldots, \bm{y}_{k+1})\in S_k(r):\bm{x}_i\neq\bm{x}_j, \bm{y}_i\neq\bm{y}_j, \forall 1\leq i<j\leq k+1\}.
$$
It is easy to see that
$$
S_k(r)=\left(\bigcup_{2\leq i<j\leq k+1}(A_{ij}\cup A'_{ij})\right)\bigcup B
$$
and
$$
|B|=|S_k(r)|-\left|\bigcup_{2\leq i<j\leq k+1}(A_{ij}\cup A'_{ij})\right|.
$$
Note that
$$
\left|\bigcup_{2\leq i<j\leq k+1}(A_{ij}\cup A'_{ij})\right|\leq\sum_{2\leq i<j\leq k+1}\left|A_{ij}\cup A'_{ij}\right|={k\choose2}\left|A_{k, k+1}\cup A'_{k, k+1}\right|
$$
and
$$
\left|A_{k, k+1}\cup A'_{k, k+1}\right|=|A_{k, k+1}|+|A'_{k, k+1}|-\left|A_{k, k+1}\cap A'_{k, k+1}\right|.
$$
We calculate
\begin{equation}
\begin{split}
|A_{k, k+1}|=&|\{(\bm{x}_1, \ldots, \bm{x}_k, \bm{x}_{k+1}, \bm{y}_1, \ldots, \bm{y}_k, \bm{y}_{k+1})\in\mathcal{E}^{2k+2}:\bm{x}_k=\bm{x}_{k+1}, \\
&\|\bm{y}_i-\bm{y}_1\|=r\|\bm{x}_i-\bm{x}_1\|\neq0, i=2, 3, \ldots, k+1\}|\\
=&|\{(\bm{x}_1, \ldots, \bm{x}_k, \bm{y}_1, \ldots, \bm{y}_k, \bm{y}_{k+1})\in\mathcal{E}^{2k+1}:\\
&\|\bm{y}_i-\bm{y}_1\|=r\|\bm{x}_i-\bm{x}_1\|\neq0, i=2, 3, \ldots, k, \text{and }\|\bm{y}_{k+1}-\bm{y}_1\|=r\|\bm{x}_k-\bm{x}_1\|\}|\\
\leq&|\{(\bm{x}_1, \ldots, \bm{x}_k, \bm{y}_1, \ldots, \bm{y}_k)\in\mathcal{E}^{2k}:\|\bm{y}_i-\bm{y}_1\|=r\|\bm{x}_i-\bm{x}_1\|\neq0, i=2, 3, \ldots, k\}|\\
&\cdot|\{\bm{y}_{k+1}:\bm{y}_{k+1}\in\mathcal{E}\}|\\
=&|\mathcal{E}|\cdot|S_{k-1}(r)|.
\end{split}
\end{equation}
Similarly, $|A'_{k, k+1}|\leq|\mathcal{E}|\cdot|S_{k-1}(r)|$ and $\left|A_{k, k+1}\cap A'_{k, k+1}\right|=|S_{k-1}(r)|$. Therefore, we have
\begin{equation}\label{bdedaxiao}
\begin{split}
|B|=&|S_k(r)|-\left|\bigcup_{2\leq i<j\leq k+1}(A_{ij}\cup A'_{ij})\right|\\
\geq&|S_k(r)|-{k\choose2}\left|A_{k, k+1}\cup A'_{k, k+1}\right|\\
\geq&|S_k(r)|-{k\choose2}\left(2|\mathcal{E}|\cdot|S_{k-1}(r)|-|S_{k-1}(r)|\right)\\
=&|S_k(r)|-{k\choose2}\left(2|\mathcal{E}|-1\right)|S_{k-1}(r)|.
\end{split}
\end{equation}

We construct an auxiliary graph $\mathcal{G}$, where $V(\mathcal{G})=\mathcal{E}\times\mathcal{E}=\{(\bm{x}, \bm{y}):\bm{x}, \bm{y}\in\mathcal{E}\}$. Two vertices $(\bm{x}, \bm{y})$ and $(\bm{x}', \bm{y}')$ in $V(\mathcal{G})$ are adjacent if and only if $\|\bm{y}'-\bm{y}\|=r\|\bm{x}'-\bm{x}\|\neq0$. It is easy to see that $\mathcal{G}$ is well defined and $|E(\mathcal{G})|=|S_1(r)|/2$. Moreover, we have
\begin{equation}
\begin{split}
|S_k(r)|=&|\{(\bm{x}_1, \bm{x}_2, \ldots, \bm{x}_{k+1}, \bm{y}_1, \bm{y}_2, \ldots, \bm{y}_{k+1})\in\mathcal{E}^{2k+2}:\|\bm{y}_i-\bm{y}_1\|=r\|\bm{x}_i-\bm{x}_1\|\neq0, i=2, 3, \ldots, k+1\}|\\
=&|\{(\bm{x}_1, \bm{y}_1, \bm{x}_2, \bm{y}_2, \ldots, \bm{x}_{k+1}, \bm{y}_{k+1})\in V(\mathcal{G})^{k+1}:(\bm{x}_1, \bm{y}_1)\text{ is adjacent to }(\bm{x}_i, \bm{y}_i), i=2, 3, \ldots, k+1\}|\\
=&\sum_{(\bm{x}_1, \bm{y}_1)\in V(\mathcal{G})}(\deg(\bm{x}_1, \bm{y}_1))^k.
\end{split}
\end{equation}
Recall the H\"{o}lder's inequality $(\sum_{i=1}^n a_i^\alpha)^{1/\alpha}(\sum_{i=1}^n b_i^\beta)^{1/\beta}\geq\sum_{i=1}^n a_ib_i$. Applying the H\"{o}lder's inequality with $a_i=(\deg(\bm{x}, \bm{y}))^{k-1}, b_i=1, n=|V(\mathcal{G})|$ and $\alpha=k/(k-1)$, it follows that
$$
|S_k(r)|^{(k-1)/k}|V(\mathcal{G})|^{1/k}\geq|S_{k-1}(r)|,
$$
in other words,
$$
|S_k(r)|\geq\frac{1}{|V(\mathcal{G})|^{1/(k-1)}}|S_{k-1}(r)|^{k/(k-1)}=\frac{1}{|\mathcal{E}|^{2/(k-1)}}|S_{k-1}(r)|^{k/(k-1)}.
$$
Applying the H\"{o}lder's inequality with $a_i=\deg(\bm{x}, \bm{y}), b_i=1, n=|V(\mathcal{G})|$ and $\alpha=k-1$, we have
$$
|S_{k-1}(r)|^{1/(k-1)}|V(\mathcal{G})|^{(k-2)/(k-1)}\geq|S_{1}(r)|,
$$
in other words,
$$
|S_{k-1}(r)|^{1/(k-1)}\geq|\mathcal{E}|^{(4-2k)/(k-1)}|S_{1}(r)|.
$$
So inequality (\ref{bdedaxiao}) becomes
\begin{equation}\label{Bxiajie}
\begin{split}
|B|\geq&|S_k(r)|-{k\choose2}\left(2|\mathcal{E}|-1\right)|S_{k-1}(r)|\\
\geq&|S_{k-1}(r)|\cdot\left(\frac{1}{|\mathcal{E}|^{2/(k-1)}}|S_{k-1}(r)|^{1/(k-1)}-2{k\choose2}|\mathcal{E}|+{k\choose2}\right)\\
\geq&|\mathcal{E}|^{-2}|S_{k-1}(r)|\cdot\left(|S_{1}(r)|-2{k\choose2}|\mathcal{E}|^3+{k\choose2}|\mathcal{E}|^2\right).
\end{split}
\end{equation}
If $|S_{1}(r)|-2{k\choose2}|\mathcal{E}|^3\geq0$, i.e. $|S_{1}(r)|\geq2{k\choose2}|\mathcal{E}|^3$, then $|S_{k-1}(r)|>0$ (since $|S_{1}(r)|>0$) and $|B|\geq{k\choose2}|S_{k-1}(r)|>0$, which means that $\mathcal{E}$ contains a pair of $k$-stars with dilation ratio $r$.

Before completing the proof of Theorem \ref{diyidingli}, we need the following theorem, which follows from (2-7), (2-9), (3-2) and (3-3) in \cite{MR3959878}.

\begin{theorem}[\cite{MR3959878}]\label{s1xiajie}
  Let $\mathcal{E}\subseteq\mathbb{F}_q^d$ with $|\mathcal{E}|\geq q^{d/2}$ and $S_1(r)$ be defined as above.
  \begin{itemize}
    \item If $d\geq2$ is even and $r\in\mathbb{F}_q^*$, then
  $$
  |S_{1}(r)|\geq q^{-1}|\mathcal{E}|^{4}-2|\mathcal{E}|^{3}-q^{d-1}|\mathcal{E}|^2-4q^{-2}|\mathcal{E}|^4-4q^{(d-2)/2}|\mathcal{E}|^{3}.
  $$
    \item If $d\geq3$ is odd and $r\in\mathbb{F}_q^+$, then
  $$
  |S_{1}(r)|\geq q^{-1}|\mathcal{E}|^{4}-2|\mathcal{E}|^{3}-2q^{d-1}|\mathcal{E}|^2-q^{-2}|\mathcal{E}|^4-2q^{(d-3)/2}|\mathcal{E}|^{3}.
  $$
  \end{itemize}
\end{theorem}

Now we are ready to complete the proof of Theorem \ref{diyidingli}.
\begin{proof}[Complete proof of Theorem \ref{diyidingli}]
By inequality (\ref{Bxiajie}), it suffices to show that $|S_{1}(r)|-2{k\choose2}|\mathcal{E}|^3\geq0$.

\textbf{Case 1}: $q\geq5$, $d\geq2$ is even and $r\in\mathbb{F}_q^*$.

Suppose $|\mathcal{E}|=tq^{d/2}$ for some $t\geq31+10{k\choose2}$.
Applying Theorem \ref{s1xiajie}, we have
\begin{equation}
\begin{split}
&|S_{1}(r)|-2{k\choose2}|\mathcal{E}|^3\\
\geq& q^{2d-1}t^4-2t^{3}q^{3d/2}-q^{2d-1}t^2-4q^{2d-2}t^4-4q^{2d-1}t^{3}-2{k\choose2}t^3q^{3d/2}\\
=&t^2q^d\left(t^2q^{d-1}-2tq^{d/2}-q^{d-1}-4q^{d-2}t^2-4q^{d-1}t-2{k\choose2}tq^{d/2}\right)\\
\geq&t^2q^d\left(t^2q^{d-1}-2tq^{d/2}-q^{d-1}-\frac{4}{5}q^{d-1}t^2-4q^{d-1}t-2{k\choose2}tq^{d/2}\right)\\
=&\frac{1}{5}t^2q^{3d/2}\left(\left(t^2-20t-5\right)q^{(d-2)/2}-10t-10{k\choose2}t\right)\\
\geq&\frac{1}{5}t^2q^{3d/2}\left(\left(t^2-21t\right)q^{(2-2)/2}-10t-10{k\choose2}t\right)\\
\geq&\frac{1}{5}t^2q^{3d/2}\left(t^2-31t-10{k\choose2}t\right)\geq0.
\end{split}
\end{equation}

\textbf{Case 2}: $d\geq3$ is odd and $r\in\mathbb{F}_q^+$.

Suppose $|\mathcal{E}|=tq^{d/2}$ for some $t\geq4+\sqrt3{k\choose2}$.
Applying Theorem \ref{s1xiajie}, we have
\begin{equation}
\begin{split}
&|S_{1}(r)|-2{k\choose2}|\mathcal{E}|^3\\
\geq& q^{2d-1}t^{4}-2t^{3}q^{3d/2}-2q^{2d-1}t^2-q^{2d-2}t^4-2q^{(4d-3)/2}t^{3}-2{k\choose2}t^3q^{3d/2}\\
=&t^2q^{2d-1}\left(t^2-2tq^{(2-d)/2}-2-q^{-1}t^2-2q^{-1/2}t-2{k\choose2}tq^{(2-d)/2}\right)\\
\geq&t^2q^{2d-1}\left(t^2-2t\cdot3^{(2-3)/2}-2-3^{-1}t^2-2\cdot3^{-1/2}t-2{k\choose2}t\cdot3^{(2-3)/2}\right)\\
=&\frac{2}{3}t^2q^{2d-1}\left(t^2-\sqrt3t-3          -\sqrt3t        -\sqrt3{k\choose2}t\right)\\
=&\frac{2}{3}t^2q^{2d-1}\left(\left(t-2\sqrt3-\sqrt3{k\choose2}\right)t-3\right)\\
\geq&\frac{2}{3}t^2q^{2d-1}\left(\left(4-2\sqrt3\right)\left(4+\sqrt3\right)-3\right)>0.
\end{split}
\end{equation}
\end{proof}

\section{Proof of Theorem \ref{dierdingli}}
In this section, we investigate the minimum size of $\mathcal{E}$ to ensure that $\mathcal{E}$ contains a pair of $4$-paths with dilation ratio $r$ and prove Theorem \ref{dierdingli}.

Let
$$
P_k(r)=\{(\bm{x}_1, \ldots, \bm{x}_{k+1}, \bm{y}_1, \ldots, \bm{y}_{k+1})\in\mathcal{E}^{2k+2}:\|\bm{y}_{i+1}-\bm{y}_i\|=r\|\bm{x}_{i+1}-\bm{x}_i\|\neq0, i=1, 2, \ldots, k\}.
$$
We consider the case $k=4$.
Let
$$
A_{1}=\{(\bm{x}_1, \bm{x}_2, \ldots, \bm{x}_{5}, \bm{y}_1, \bm{y}_2, \ldots, \bm{y}_{5})\in P_4(r):\bm{x}_1=\bm{x}_3\},
$$
$$
A'_{1}=\{(\bm{x}_1, \bm{x}_2, \ldots, \bm{x}_{5}, \bm{y}_1, \bm{y}_2, \ldots, \bm{y}_{5})\in P_4(r):\bm{y}_1=\bm{y}_3\},
$$
$$
A_{2}=\{(\bm{x}_1, \bm{x}_2, \ldots, \bm{x}_{5}, \bm{y}_1, \bm{y}_2, \ldots, \bm{y}_{5})\in P_4(r):\bm{x}_1=\bm{x}_4\},
$$
$$
A'_{2}=\{(\bm{x}_1, \bm{x}_2, \ldots, \bm{x}_{5}, \bm{y}_1, \bm{y}_2, \ldots, \bm{y}_{5})\in P_4(r):\bm{y}_1=\bm{y}_4\},
$$
$$
A_{3}=\{(\bm{x}_1, \bm{x}_2, \ldots, \bm{x}_{5}, \bm{y}_1, \bm{y}_2, \ldots, \bm{y}_{5})\in P_4(r):\bm{x}_1=\bm{x}_5\},
$$
$$
A'_{3}=\{(\bm{x}_1, \bm{x}_2, \ldots, \bm{x}_{5}, \bm{y}_1, \bm{y}_2, \ldots, \bm{y}_{5})\in P_4(r):\bm{y}_1=\bm{y}_5\},
$$
$$
A_{4}=\{(\bm{x}_1, \bm{x}_2, \ldots, \bm{x}_{5}, \bm{y}_1, \bm{y}_2, \ldots, \bm{y}_{5})\in P_4(r):\bm{x}_3=\bm{x}_5\},
$$
$$
A'_{4}=\{(\bm{x}_1, \bm{x}_2, \ldots, \bm{x}_{5}, \bm{y}_1, \bm{y}_2, \ldots, \bm{y}_{5})\in P_4(r):\bm{y}_3=\bm{y}_5\},
$$
$$
A_{5}=\{(\bm{x}_1, \bm{x}_2, \ldots, \bm{x}_{5}, \bm{y}_1, \bm{y}_2, \ldots, \bm{y}_{5})\in P_4(r):\bm{x}_2=\bm{x}_5\},
$$
$$
A'_{5}=\{(\bm{x}_1, \bm{x}_2, \ldots, \bm{x}_{5}, \bm{y}_1, \bm{y}_2, \ldots, \bm{y}_{5})\in P_4(r):\bm{y}_2=\bm{y}_5\},
$$
$$
A_{6}=\{(\bm{x}_1, \bm{x}_2, \ldots, \bm{x}_{5}, \bm{y}_1, \bm{y}_2, \ldots, \bm{y}_{5})\in P_4(r):\bm{x}_2=\bm{x}_4\},
$$
$$
A'_{6}=\{(\bm{x}_1, \bm{x}_2, \ldots, \bm{x}_{5}, \bm{y}_1, \bm{y}_2, \ldots, \bm{y}_{5})\in P_4(r):\bm{y}_2=\bm{y}_4\},
$$
and
$$
C=\{(\bm{x}_1, \bm{x}_2, \ldots, \bm{x}_{5}, \bm{y}_1, \bm{y}_2, \ldots, \bm{y}_{5})\in P_4(r):\bm{x}_i\neq\bm{x}_j, \bm{y}_i\neq\bm{y}_j, \forall 1\leq i<j\leq 5\}.
$$
It is easy to see that
$$
P_4(r)=\left(\bigcup_{i=1}^6(A_{i}\cup A'_{i})\right)\bigcup C
$$
and
$$
|C|=|P_4(r)|-\left|\bigcup_{i=1}^6(A_{i}\cup A'_{i})\right|\geq|P_4(r)|-\sum_{i=1}^6\left|A_{i}\cup A'_{i}\right|.
$$

We calculate
\begin{equation}
\begin{split}
|A_{1}|=&|\{(\bm{x}_1, \ldots, \bm{x}_{5}, \bm{y}_1, \ldots, \bm{y}_{5})\in\mathcal{E}^{10}:\|\bm{y}_{i+1}-\bm{y}_i\|=r\|\bm{x}_{i+1}-\bm{x}_i\|\neq0, i=1, 2, 3, 4, \bm{x}_1=\bm{x}_{3}\}|\\
=&|\{(\bm{x}_2, \ldots, \bm{x}_5, \bm{y}_1, \ldots, \bm{y}_{5})\in\mathcal{E}^{9}:\\
&\|\bm{y}_{i+1}-\bm{y}_i\|=r\|\bm{x}_{i+1}-\bm{x}_i\|\neq0, i=2, 3, 4, \text{and }\|\bm{y}_{2}-\bm{y}_1\|=r\|\bm{x}_2-\bm{x}_3\|\}|\\
\leq&|\{(\bm{x}_2, \ldots, \bm{x}_5, \bm{y}_2, \ldots, \bm{y}_5)\in\mathcal{E}^{8}:\|\bm{y}_{i+1}-\bm{y}_i\|=r\|\bm{x}_{i+1}-\bm{x}_i\|\neq0, i=2, 3, 4\}|\\
&\cdot|\{\bm{y}_{1}:\bm{y}_{1}\in\mathcal{E}\}|\\
=&|\mathcal{E}|\cdot|P_3(r)|.
\end{split}
\end{equation}
Similarly, $|A_{i}|, |A'_i|\leq|\mathcal{E}|\cdot|P_3(r)|$ for $i=1, 2, \ldots, 5$.
On the other hand, we have
\begin{equation}\label{a4binga4}
\begin{split}
|A_{6}\cup A'_6|=&|\{(\bm{x}_1, \ldots, \bm{x}_{5}, \bm{y}_1, \ldots, \bm{y}_{5})\in P_4(r):\bm{x}_2=\bm{x}_4\text{ or }\bm{y}_2=\bm{y}_4\}|\\
=&|\{(\bm{x}_1, \ldots, \bm{x}_{5}, \bm{y}_1, \ldots, \bm{y}_{5})\in P_4(r):\bm{x}_2=\bm{x}_4\text{ and }\bm{y}_2\neq\bm{y}_4\}|\\
&+|\{(\bm{x}_1, \ldots, \bm{x}_{5}, \bm{y}_1, \ldots, \bm{y}_{5})\in P_4(r):\bm{x}_2\neq\bm{x}_4\text{ and }\bm{y}_2=\bm{y}_4\}|\\
&+|\{(\bm{x}_1, \ldots, \bm{x}_{5}, \bm{y}_1, \ldots, \bm{y}_{5})\in P_4(r):\bm{x}_2=\bm{x}_4\text{ and }\bm{y}_2=\bm{y}_4\}|\\
:=&I+II+III.
\end{split}
\end{equation}
We calculate
\begin{equation}
\begin{split}
I=&|\{(\bm{x}_1, \ldots, \bm{x}_{5}, \bm{y}_1, \ldots, \bm{y}_{5})\in\mathcal{E}^{10}:\|\bm{y}_{i+1}-\bm{y}_i\|=r\|\bm{x}_{i+1}-\bm{x}_i\|\neq0, i=1, 2, 3, 4, \bm{x}_2=\bm{x}_{4}, \bm{y}_2\neq\bm{y}_4\}|\\
=&|\{(\bm{x}_1, \ldots, \bm{x}_5, \bm{y}_1, \ldots, \bm{y}_{5})\in\mathcal{E}^{10}:\bm{x}_2=\bm{x}_{4}, \bm{y}_2\neq\bm{y}_4, \\
&\|\bm{y}_{i+1}-\bm{y}_i\|=r\|\bm{x}_{i+1}-\bm{x}_i\|\neq0, i=1, 2, 4,\|\bm{y}_{4}-\bm{y}_3\|=r\|\bm{x}_2-\bm{x}_3\|\}|\\
\leq&|\{(\bm{x}_1, \bm{x}_2, \bm{x}_4, \bm{x}_5, \bm{y}_1, \ldots, \bm{y}_5)\in\mathcal{E}^{9}:\bm{x}_2=\bm{x}_{4}, \bm{y}_{2}\neq\bm{y}_4, \\
&\|\bm{y}_{2}-\bm{y}_1\|=r\|\bm{x}_{2}-\bm{x}_1\|\neq0, \|\bm{y}_{3}-\bm{y}_2\|=\|\bm{y}_{4}-\bm{y}_3\|\neq0, \|\bm{y}_{5}-\bm{y}_4\|=r\|\bm{x}_{5}-\bm{x}_4\|\neq0\}|\\
&\cdot|\{\bm{x}_{3}:\bm{x}_{3}\in\mathcal{E}\}|\\
=&|\{(\bm{x}_1, \bm{x}_2, \bm{x}_4, \bm{x}_5, \bm{y}_1, \ldots, \bm{y}_5)\in\mathcal{E}^{9}:\bm{x}_2=\bm{x}_{4}, \bm{y}_{2}\neq\bm{y}_4, \\
&\|\bm{y}_{2}-\bm{y}_1\|=r\|\bm{x}_{2}-\bm{x}_1\|\neq0, \|\bm{y}_{3}-\bm{y}_2\|=\|\bm{y}_{4}-\bm{y}_3\|\neq0, \|\bm{y}_{5}-\bm{y}_4\|=r\|\bm{x}_{5}-\bm{x}_4\|\neq0\}|\\
&\cdot|\mathcal{E}|.
\end{split}
\end{equation}
If we write $\bm{y}_2=(u_1, u_2, \ldots, u_d), \bm{y}_3=(v_1, v_2, \ldots, v_d), \bm{y}_4=(w_1, w_2, \ldots, w_d)$, where $(u_1, u_2, \ldots, u_d)\neq(w_1, w_2, \ldots, w_d)$, then $\|\bm{y}_{3}-\bm{y}_2\|=\|\bm{y}_{4}-\bm{y}_3\|$ becomes
$$
(v_1-u_1)^2+(v_2-u_2)^2+\cdots+(v_d-u_d)^2=(v_1-w_1)^2+(v_2-w_2)^2+\cdots+(v_d-w_d)^2,
$$i.e.
\begin{equation}\label{zhixian}
2(w_1-u_1)v_1+2(w_2-u_2)v_2+\cdots+2(w_d-u_d)v_d=w_1^2+w_2^2+\cdots+w_d^2-u_1^2-u_2^2-\cdots-u_d^2.
\end{equation}
Note that $\bm{y}_{2}\neq\bm{y}_4$. So there are at most $q^{d-1}$ solutions to Equation (\ref{zhixian}) for $\bm{y}_3$ whenever $\bm{y}_2\neq\bm{y}_4$. Recall that $\bm{y}_3\in \mathcal{E}$. Thus
\begin{equation}
\begin{split}
&|\{(\bm{x}_1, \bm{x}_2, \bm{x}_4, \bm{x}_5, \bm{y}_1, \ldots, \bm{y}_5)\in\mathcal{E}^{9}:\bm{x}_2=\bm{x}_{4}, \bm{y}_{2}\neq\bm{y}_4,\\
&\|\bm{y}_{2}-\bm{y}_1\|=r\|\bm{x}_{2}-\bm{x}_1\|\neq0, \|\bm{y}_{3}-\bm{y}_2\|=\|\bm{y}_{4}-\bm{y}_3\|\neq0, \|\bm{y}_{5}-\bm{y}_4\|=r\|\bm{x}_{5}-\bm{x}_4\|\neq0\}|\\
\leq&|\{(\bm{x}_1, \bm{x}_2, \bm{x}_4, \bm{x}_5, \bm{y}_1, \bm{y}_2, \bm{y}_4, \bm{y}_5)\in\mathcal{E}^{8}:\bm{x}_2=\bm{x}_{4}, \bm{y}_{2}\neq\bm{y}_4, \\ &\|\bm{y}_{2}-\bm{y}_1\|=r\|\bm{x}_{2}-\bm{x}_1\|\neq0, \|\bm{y}_{5}-\bm{y}_4\|=r\|\bm{x}_{5}-\bm{x}_4\|\neq0\}|\\
&\cdot \min\left\{q^{d-1}, |\mathcal{E}|\right\}\\
=&|\{(\bm{x}_1, \bm{x}_2, \bm{x}_5, \bm{y}_1, \bm{y}_2, \bm{y}_4, \bm{y}_5)\in\mathcal{E}^{7}:\bm{y}_{2}\neq\bm{y}_4, \|\bm{y}_{2}-\bm{y}_1\|=r\|\bm{x}_{2}-\bm{x}_1\|\neq0, \|\bm{y}_{5}-\bm{y}_4\|=r\|\bm{x}_{5}-\bm{x}_2\|\neq0\}|\\
&\cdot \min\left\{q^{d-1}, |\mathcal{E}|\right\}\\
\leq&|\{(\bm{x}_1, \bm{x}_2, \bm{x}_5, \bm{y}_1, \bm{y}_2, \bm{y}_4, \bm{y}_5)\in\mathcal{E}^{7}:\|\bm{y}_{2}-\bm{y}_1\|=r\|\bm{x}_{2}-\bm{x}_1\|\neq0, \|\bm{y}_{5}-\bm{y}_4\|=r\|\bm{x}_{5}-\bm{x}_2\|\neq0\}|\\
&\cdot \min\left\{q^{d-1}, |\mathcal{E}|\right\}.
\end{split}
\end{equation}
In order to obtain an upper bound for the size of the set
$$
\{(\bm{x}_1, \bm{x}_2, \bm{x}_5, \bm{y}_1, \bm{y}_2, \bm{y}_4, \bm{y}_5)\in\mathcal{E}^{7}:\|\bm{y}_{2}-\bm{y}_1\|=r\|\bm{x}_{2}-\bm{x}_1\|\neq0, \|\bm{y}_{5}-\bm{y}_4\|=r\|\bm{x}_{5}-\bm{x}_2\|\neq0\},
$$
we first choose $(\bm{x}_1, \bm{x}_2, \bm{y}_1, \bm{y}_2)\in\mathcal{E}^4$, such that $\|\bm{y}_{2}-\bm{y}_1\|=r\|\bm{x}_{2}-\bm{x}_1\|\neq0$. So there are $|P_1(r)|$ choices for $(\bm{x}_1, \bm{x}_2, \bm{y}_1, \bm{y}_2)$. Then we choose $\bm{x}_5$ and $\bm{y}_4$ from $\mathcal{E}$ such that $\|\bm{x}_5-\bm{x}_2\|\neq0$ and there are at most $|\mathcal{E}|^2$ choices. Finally, we choose $\bm{y}_5\in\mathcal{E}$ such that $\|\bm{y}_{5}-\bm{y}_4\|=r\|\bm{x}_{5}-\bm{x}_2\|$. If we use $S(\bm{x}, t)$ to denote the sphere with center $\bm{x}$ and radius $t$, i.e. $S(\bm{x}, t):=\{\bm{y}\in\mathbb{F}_q^d: \|\bm{y}-\bm{x}\|=t\}$, then $\bm{y}_{5}\in S(\bm{y}_4, r\|\bm{x}_{5}-\bm{x}_2\|)$. The following theorem gives the exact value of $|S(\bm{x}, t)|$, which can be deduced by Theorems 6.26 and 6.27 in \cite{MR1429394}.
\begin{theorem}[\cite{MR1429394}]\label{qiumian}
  Let $S(\bm{x}, t)$ be the sphere with center $\bm{x}$ and radius $t$ in $\mathbb{F}_q^d$.
  \begin{itemize}
    \item If $d$ is even, then
  $$
  |S(\bm{x}, t)|=q^{d-1}+\mu(t)q^{(d-2)/2}\psi\left((-1)^{d/2}\right),
  $$
  where $\mu(t)=q-1$ if $t=0$ and $\mu(t)=-1$ if $t\in\mathbb{F}_q^*$, and $\psi$ is a quadratic character of $\mathbb{F}_q$.
    \item If $d\geq3$ is odd, then
  $$
  |S(\bm{x}, t)|=q^{d-1}+q^{(d-1)/2}\eta\left((-1)^{(d-1)/2}t\right),
  $$
  where $\eta$ is a quadratic character of $\mathbb{F}_q^*$ and $\eta(0)=0$.
  \end{itemize}
\end{theorem}
If $d=2$, applying Theorem \ref{qiumian} with $\bm{x}=\bm{y}_4$ and $t=r\|\bm{x}_{5}-\bm{x}_2\|\neq0$, we obtain that there are at most $q^{d-1}+q^{(d-2)/2}$ choices for $\bm{y}_5$, where $q^{d-1}+q^{(d-2)/2}=q^{d-1}\left(1+q^{-d/2}\right)\leq q^{d-1}\left(1+\frac{1}{3}\right)\leq1.5q^{d-1}$. If $d\geq3$, applying Theorem \ref{qiumian} with $\bm{x}=\bm{y}_4$ and $t=r\|\bm{x}_{5}-\bm{x}_2\|\neq0$, we obtain that there are at most $q^{d-1}+q^{(d-1)/2}$ choices for $\bm{y}_5$, where $q^{d-1}+q^{(d-1)/2}=q^{d-1}\left(1+q^{(1-d)/2}\right)\leq q^{d-1}\left(1+\frac{1}{3}\right)\leq1.5q^{d-1}$. So there are always at most $1.5q^{d-1}$ choices for $\bm{y}_5$. On the other hand, $\bm{y}_5\in\mathcal{E}$ implies that there are at most $|\mathcal{E}|$ choices for $\bm{y}_5$.
Therefore, we have
\begin{equation}
  I\leq|\mathcal{E}|\cdot \min\left\{q^{d-1}, |\mathcal{E}|\right\}\cdot|P_1(r)|\cdot|\mathcal{E}|^2\cdot\min\left\{1.5q^{d-1}, |\mathcal{E}|\right\}=|\mathcal{E}|^3\cdot|P_1(r)|\cdot \min\left\{q^{d-1}, |\mathcal{E}|\right\}\cdot\min\left\{1.5q^{d-1}, |\mathcal{E}|\right\}.
\end{equation}
When $2\leq d\leq4$, $\min\left\{q^{d-1}, |\mathcal{E}|\right\}\cdot\min\left\{1.5q^{d-1}, |\mathcal{E}|\right\}\leq q^{d-1}\cdot1.5q^{d-1}=1.5q^{2d-2}$. When $d\geq5$, $\min\left\{q^{d-1}, |\mathcal{E}|\right\}\cdot\min\left\{1.5q^{d-1}, |\mathcal{E}|\right\}\leq|\mathcal{E}|\cdot|\mathcal{E}|=|\mathcal{E}|^2$. So
\begin{equation}\label{dange}
  I\leq\left\{
\begin{array}{lllll}
  1.5|\mathcal{E}|^3\cdot q^{2d-2}\cdot|P_1(r)|, & \text{ if }2\leq d\leq4;\\
  |\mathcal{E}|^5\cdot|P_1(r)|, & \text{ if }d\geq5.
\end{array}
  \right.
\end{equation}
Similarly,
\begin{equation}\label{dange2}
  II\leq\left\{
\begin{array}{lllll}
  1.5|\mathcal{E}|^3\cdot q^{2d-2}\cdot|P_1(r)|, & \text{ if }2\leq d\leq4;\\
  |\mathcal{E}|^5\cdot|P_1(r)|, & \text{ if }d\geq5.
\end{array}
  \right.
\end{equation}
Furthermore, we calculate
\begin{equation}\label{dange3}
\begin{split}
III=&|\{(\bm{x}_1, \bm{x}_2, \bm{x}_3, \bm{x}_5, \bm{y}_1, \bm{y}_2, \bm{y}_3, \bm{y}_5)\in\mathcal{E}^{8}:\\
&\|\bm{y}_{2}-\bm{y}_1\|=r\|\bm{x}_{2}-\bm{x}_1\|\neq0, \|\bm{y}_{3}-\bm{y}_2\|=r\|\bm{x}_{3}-\bm{x}_2\|\neq0, \|\bm{y}_{5}-\bm{y}_2\|=r\|\bm{x}_{5}-\bm{x}_2\|\neq0\}|\\
\leq&|\{(\bm{x}_1, \bm{x}_2, \bm{x}_5, \bm{y}_1, \bm{y}_2, \bm{y}_5)\in\mathcal{E}^{6}:\|\bm{y}_{2}-\bm{y}_1\|=r\|\bm{x}_{2}-\bm{x}_1\|\neq0, \|\bm{y}_{5}-\bm{y}_2\|=r\|\bm{x}_{5}-\bm{x}_2\|\neq0\}|\\
&\cdot|\{\bm{x}_{3}:\bm{x}_{3}\in\mathcal{E}\}|\cdot|\{\bm{y}_{3}:\bm{y}_{3}\in\mathcal{E}\}|\\
=&|\mathcal{E}|^2\cdot|P_2(r)|.
\end{split}
\end{equation}
Combining equations (\ref{a4binga4}), (\ref{dange}), (\ref{dange2}), and (\ref{dange3}), we conclude that
\begin{equation}
  |A_{6}\cup A'_6|\leq\left\{
\begin{array}{lllll}
  3|\mathcal{E}|^3\cdot q^{2d-2}\cdot|P_1(r)|+|\mathcal{E}|^2\cdot|P_2(r)|, & \text{ if }2\leq d\leq4;\\
  2|\mathcal{E}|^5\cdot|P_1(r)|+|\mathcal{E}|^2\cdot|P_2(r)|, & \text{ if }d\geq5.
\end{array}
  \right.
\end{equation}

Hence,
\begin{equation}\label{Cxiajie}
\begin{split}
|C|\geq&|P_4(r)|-\sum_{i=1}^6\left|A_{i}\cup A'_{i}\right|\\
\geq&|P_4(r)|-\sum_{i=1}^5\left(|A_{i}|+|A'_{i}|\right)-\left|A_{6}\cup A'_{6}\right|\\
\geq&\left\{
\begin{array}{lllll}
  |P_4(r)|-10|\mathcal{E}|\cdot|P_3(r)|-3|\mathcal{E}|^3\cdot q^{2d-2}\cdot|P_1(r)|-|\mathcal{E}|^2\cdot|P_2(r)|, & \text{ if }2\leq d\leq4;\\
  |P_4(r)|-10|\mathcal{E}|\cdot|P_3(r)|-2|\mathcal{E}|^5\cdot|P_1(r)|-|\mathcal{E}|^2\cdot|P_2(r)|, & \text{ if }d\geq5.
\end{array}
  \right.
\end{split}
\end{equation}

We construct an auxiliary graph $\mathcal{G}$, where $V(\mathcal{G})=\mathcal{E}\times\mathcal{E}=\{(\bm{x}, \bm{y}):\bm{x}, \bm{y}\in\mathcal{E}\}$. Two vertices $(\bm{x}, \bm{y})$ and $(\bm{x}', \bm{y}')$ are adjacent if and only if $\|\bm{y}'-\bm{y}\|=r\|\bm{x}'-\bm{x}\|\neq0$. We have
\begin{equation}
\begin{split}
|P_k(r)|=&|\{(\bm{x}_1, \ldots, \bm{x}_{k+1}, \bm{y}_1, \ldots, \bm{y}_{k+1})\in\mathcal{E}^{2k+2}:\|\bm{y}_{i+1}-\bm{y}_i\|=r\|\bm{x}_{i+1}-\bm{x}_i\|\neq0, 1\leq i\leq k\}|\\
=&|\{(\bm{x}_1, \bm{y}_1, \bm{x}_2, \bm{y}_2, \ldots, \bm{x}_{k+1}, \bm{y}_{k+1})\in V(\mathcal{G})^{k+1}:(\bm{x}_i, \bm{y}_i)\text{ is adjacent to }(\bm{x}_{i+1}, \bm{y}_{i+1}), 1\leq i\leq k\}|\\
:=&w_k(\mathcal{G}),
\end{split}
\end{equation}
where we use $w_k(\mathcal{G})$ to denote the number of  walks of length $k$ in $\mathcal{G}$. The following theorems give the relations between $w_k(\mathcal{G})$ for different $k$, which were proved in \cite{MR184950} and \cite{lagarias1983inequality}, respectively.
\begin{theorem}[\cite{MR184950}]
  For a graph $G$  and a positive integer $k$, we have
  \begin{equation}\label{walkbudengshi2}
      w_{1}(G)^k\leq w_{0}(G)^{k-1}\cdot w_{k}(G).
  \end{equation}
\end{theorem}
\begin{theorem}[\cite{lagarias1983inequality}]
  For a graph $G$  and nonnegative integers $a, b$, we have
  \begin{equation}\label{walkbudengshi}
      w_{2a+b}(G)\cdot w_b(G)\leq w_{0}(G)\cdot w_{2(a+b)}(G).
  \end{equation}
\end{theorem}
Applying inequality (\ref{walkbudengshi2}) with $k=4$, we have $|P_1(r)|^4\leq|P_0(r)|^3\cdot|P_4(r)|$. Applying inequality (\ref{walkbudengshi2}) with $k=2$, we have $|P_1(r)|^2\leq|P_0(r)|\cdot|P_2(r)|$. Applying inequality (\ref{walkbudengshi}) with $a=b=1$, we have $|P_3(r)|\cdot|P_1(r)|\leq|P_0(r)|\cdot|P_4(r)|$. Applying inequality (\ref{walkbudengshi}) with $a=0$ and $b=2$, we have $|P_2(r)|^2\leq|P_0(r)|\cdot|P_4(r)|$. Note that $|P_0(r)|$ is exactly the number of vertices in $\mathcal{G}$. So
\begin{equation}\label{p4p1}
  |P_4(r)|\geq|\mathcal{E}|^{-6}|P_1(r)|^4,
\end{equation}
\begin{equation}\label{p2p1}
  |P_2(r)|\geq|\mathcal{E}|^{-2}|P_1(r)|^2,
\end{equation}
\begin{equation}\label{p4p3}
  |P_4(r)|\geq|\mathcal{E}|^{-2}|P_3(r)|\cdot|P_1(r)|,
\end{equation}
and
\begin{equation}\label{p4p2}
  |P_4(r)|\geq|\mathcal{E}|^{-2}|P_2(r)|^2.
\end{equation}
\subsection{The case $2\leq d\leq4$}\label{24zhijian}
In the case $2\leq d\leq4$, inequality (\ref{Cxiajie}) becomes
\begin{equation}
\begin{split}
|C|\geq&|P_4(r)|-10|\mathcal{E}|\cdot|P_3(r)|-3|\mathcal{E}|^3\cdot q^{2d-2}\cdot|P_1(r)|-|\mathcal{E}|^2\cdot|P_2(r)|\\
\geq&\left(0.96536|P_4(r)|-10|\mathcal{E}|\cdot|P_3(r)|\right)+\left(0.02357|P_4(r)|-3|\mathcal{E}|^3\cdot q^{2d-2}\cdot|P_1(r)|\right)\\
&+\left(0.01105|P_4(r)|-|\mathcal{E}|^2\cdot|P_2(r)|\right)\\
:=&I'+II'+III'.
\end{split}
\end{equation}

Applying inequality (\ref{p4p3}), we have
\begin{equation}\label{iyipie}
\begin{split}
I'=&0.96536|P_4(r)|-10|\mathcal{E}|\cdot|P_3(r)|\\
\geq&|\mathcal{E}|^{-2}|P_3(r)|\left(0.96536|P_1(r)|-10|\mathcal{E}|^3\right).
\end{split}
\end{equation}
Applying inequality (\ref{p4p1}), we have
\begin{equation}\label{iiyipie}
\begin{split}
II'=&0.02357|P_4(r)|-3|\mathcal{E}|^3\cdot q^{2d-2}\cdot|P_1(r)|\\
\geq&|\mathcal{E}|^{-6}|P_1(r)|\left(0.02357|P_1(r)|^3-3q^{2d-2}|\mathcal{E}|^9\right).
\end{split}
\end{equation}
Applying inequalities (\ref{p4p2}) and (\ref{p2p1}), we have
\begin{equation}\label{iiiyipie}
\begin{split}
III'=&0.01105|P_4(r)|-|\mathcal{E}|^2\cdot|P_2(r)|\\
\geq&|\mathcal{E}|^{-2}|P_2(r)|\left(0.01105|P_2(r)|-|\mathcal{E}|^4\right)\\
\geq&|\mathcal{E}|^{-4}|P_2(r)|\left(0.01105|P_1(r)|^2-|\mathcal{E}|^6\right).
\end{split}
\end{equation}

Assume that $|P_1(r)|>5.03023q^{(2d-2)/3}|\mathcal{E}|^3$. Then  $0.02357|P_1(r)|^3-3q^{2d-2}|\mathcal{E}|^9>0$, so we have $II'>0$ in inequality (\ref{iiyipie}). Moreover,
\begin{equation}
\begin{split}
&0.96536|P_1(r)|-10|\mathcal{E}|^3\\
>&0.96536\cdot5.03023q^{(2d-2)/3}|\mathcal{E}|^3-10|\mathcal{E}|^3\\
\geq&0.96536\cdot5.03023\cdot 3^{(2\cdot2-2)/3}|\mathcal{E}|^3-10|\mathcal{E}|^3\\
>&0.1|\mathcal{E}|^3>0.
\end{split}
\end{equation}
So $I'>0$ in inequality (\ref{iyipie}).
\begin{equation}
\begin{split}
&\sqrt{0.01105}|P_1(r)|-|\mathcal{E}|^3\\
>&\sqrt{0.01105}\cdot5.03023q^{(2d-2)/3}|\mathcal{E}|^3-|\mathcal{E}|^3\\
\geq&\sqrt{0.01105}\cdot5.03023\cdot 3^{(2\cdot2-2)/3}|\mathcal{E}|^3-|\mathcal{E}|^3\\
>&0.09|\mathcal{E}|^3>0.
\end{split}
\end{equation}
So $III'>0$ in inequality (\ref{iiiyipie}).

\subsection{The case $d=5$}\label{ddengyu5}
In the case $d=5$, inequality (\ref{Cxiajie}) becomes
\begin{equation}
\begin{split}
|C|\geq&|P_4(r)|-10|\mathcal{E}|\cdot|P_3(r)|-2|\mathcal{E}|^5\cdot|P_1(r)|-|\mathcal{E}|^2\cdot|P_2(r)|\\
\geq&\left(0.16|P_4(r)|-10|\mathcal{E}|\cdot|P_3(r)|\right)+\left(0.82|P_4(r)|-2|\mathcal{E}|^5\cdot|P_1(r)|\right)+\left(0.01|P_4(r)|-|\mathcal{E}|^2\cdot|P_2(r)|\right)\\
:=&I''+II''+III''.
\end{split}
\end{equation}
As in Subsection \ref{24zhijian}, we have
\begin{equation}\label{ierpie}
\begin{split}
I''=&0.16|P_4(r)|-10|\mathcal{E}|\cdot|P_3(r)|\\
\geq&|\mathcal{E}|^{-2}|P_3(r)|\left(0.16|P_1(r)|-10|\mathcal{E}|^3\right),
\end{split}
\end{equation}
\begin{equation}\label{iierpie}
\begin{split}
II''=&0.82|P_4(r)|-2|\mathcal{E}|^5\cdot|P_1(r)|\\
\geq&|\mathcal{E}|^{-6}|P_1(r)|\left(0.82|P_1(r)|^3-2|\mathcal{E}|^{11}\right),
\end{split}
\end{equation}
and
\begin{equation}\label{iiierpie}
\begin{split}
III''=&0.01|P_4(r)|-|\mathcal{E}|^2\cdot|P_2(r)|\\
\geq&|\mathcal{E}|^{-2}|P_2(r)|\left(0.01|P_2(r)|-|\mathcal{E}|^4\right)\\
\geq&|\mathcal{E}|^{-4}|P_2(r)|\left(0.01|P_1(r)|^2-|\mathcal{E}|^6\right).
\end{split}
\end{equation}

Assume that $|P_1(r)|>1.34609|\mathcal{E}|^{11/3}$ and $|\mathcal{E}|\geq12q^3$. Then $0.82|P_1(r)|^3-2|\mathcal{E}|^{11}>0$, and we have $II''>0$ in inequality (\ref{iierpie}). Moreover, $|P_1(r)|\geq1.34609\cdot\left(12q^3\right)^{2/3}|\mathcal{E}|^{3}\geq63.49|\mathcal{E}|^{3}$. So
\begin{equation}
\begin{split}
&0.16|P_1(r)|-10|\mathcal{E}|^3\\
>&0.16\cdot63.49|\mathcal{E}|^3-10|\mathcal{E}|^3\\
>&0.15|\mathcal{E}|^3>0,
\end{split}
\end{equation}
and $I''>0$ in inequality (\ref{ierpie}).
\begin{equation}
\begin{split}
&\sqrt{0.01}|P_1(r)|-|\mathcal{E}|^3\\
>&\sqrt{0.01}\cdot63.49|\mathcal{E}|^3-|\mathcal{E}|^3\\
>&5|\mathcal{E}|^3>0,
\end{split}
\end{equation}
and $III''>0$ in inequality (\ref{iiierpie}).
\subsection{The case $d\geq6$}\label{ddayudengyu6}
In the case $d\geq6$, inequality (\ref{Cxiajie}) becomes
\begin{equation}
\begin{split}
|C|\geq&|P_4(r)|-10|\mathcal{E}|\cdot|P_3(r)|-2|\mathcal{E}|^5\cdot|P_1(r)|-|\mathcal{E}|^2\cdot|P_2(r)|\\
\geq&\left(0.0191|P_4(r)|-10|\mathcal{E}|\cdot|P_3(r)|\right)+\left(0.98|P_4(r)|-2|\mathcal{E}|^5\cdot|P_1(r)|\right)+\left(0.0009|P_4(r)|-|\mathcal{E}|^2\cdot|P_2(r)|\right)\\
:=&I'''+II'''+III'''.
\end{split}
\end{equation}
As in Subsection \ref{24zhijian}, we have
\begin{equation}\label{isanpie}
\begin{split}
I'''=&0.0191|P_4(r)|-10|\mathcal{E}|\cdot|P_3(r)|\\
\geq&|\mathcal{E}|^{-2}|P_3(r)|\left(0.0191|P_1(r)|-10|\mathcal{E}|^3\right),
\end{split}
\end{equation}
\begin{equation}\label{iisanpie}
\begin{split}
II'''=&0.98|P_4(r)|-2|\mathcal{E}|^5\cdot|P_1(r)|\\
\geq&|\mathcal{E}|^{-6}|P_1(r)|\left(0.98|P_1(r)|^3-2|\mathcal{E}|^{11}\right),
\end{split}
\end{equation}
and
\begin{equation}\label{iiisanpie}
\begin{split}
III'''=&0.0009|P_4(r)|-|\mathcal{E}|^2\cdot|P_2(r)|\\
\geq&|\mathcal{E}|^{-2}|P_2(r)|\left(0.0009|P_2(r)|-|\mathcal{E}|^4\right)\\
\geq&|\mathcal{E}|^{-4}|P_2(r)|\left(0.0009|P_1(r)|^2-|\mathcal{E}|^6\right).
\end{split}
\end{equation}

Assume that $|P_1(r)|>1.26844|\mathcal{E}|^{11/3}$ and $|\mathcal{E}|\geq313q^{d/2}$. Then $0.98|P_1(r)|^3-2|\mathcal{E}|^{11}>0$, and we have $II'''>0$ in inequality (\ref{iisanpie}). Moreover, $|P_1(r)|\geq1.26844\cdot\left(313q^3\right)^{2/3}|\mathcal{E}|^{3}\geq526.27|\mathcal{E}|^{3}$. So
\begin{equation}
\begin{split}
&0.0191|P_1(r)|-10|\mathcal{E}|^3\\
>&0.0191\cdot526.27|\mathcal{E}|^3-10|\mathcal{E}|^3\\
>&0.05|\mathcal{E}|^3>0,
\end{split}
\end{equation}
and $I'''>0$ in inequality (\ref{isanpie}).
\begin{equation}
\begin{split}
&\sqrt{0.0009}|P_1(r)|-|\mathcal{E}|^3\\
>&\sqrt{0.0009}\cdot526.27|\mathcal{E}|^3-|\mathcal{E}|^3\\
>&10|\mathcal{E}|^3>0,
\end{split}
\end{equation}
and $III'''>0$ in inequality (\ref{iiisanpie}).

\subsection{Completing the proof of Theorem \ref{dierdingli}}
Now we are ready to complete the proof of Theorem \ref{dierdingli}.
\begin{proof}[Complete proof of Theorem \ref{dierdingli}]

\textbf{Case 1}: $q\geq5$, $2\leq d\leq4$ is even and $r\in\mathbb{F}_q^*$.

As in Subsection \ref{24zhijian}, it suffices to prove $|P_1(r)|>5.03023q^{(2d-2)/3}|\mathcal{E}|^3$. If so, then $|C|\geq I'+II'+III'>0$. And hence $\mathcal{E}$ contains a pair of $4$-paths with dilation ratio $r$.

Suppose $|\mathcal{E}|=tq^{(2d+1)/3}$ for some $t\geq36$. Note that $P_1(r)=S_1(r)$. So we can apply Theorem \ref{s1xiajie} and obtain that
  \begin{equation}
\begin{split}
  &|P_1(r)|-5.03023q^{(2d-2)/3}|\mathcal{E}|^3\\
  \geq&q^{-1}|\mathcal{E}|^{4}-2|\mathcal{E}|^{3}-q^{d-1}|\mathcal{E}|^2-4q^{-2}|\mathcal{E}|^4-4q^{(d-2)/2}|\mathcal{E}|^{3}-5.03023q^{(2d-2)/3}|\mathcal{E}|^3\\
  \geq&t^{4}q^{(8d+1)/3}-2t^{3}q^{2d+1}-t^2q^{(7d-1)/3}-4t^4q^{(8d-2)/3}-4t^{3}q^{5d/2}-5.03023t^3q^{(8d+1)/3}\\
  =&t^2q^{(8d+1)/3}\left(t^2-2tq^{(-2d+2)/3}-q^{(-d-2)/3}-4t^2q^{-1}    -4tq^{(-d-2)/6}-5.03023t\right)\\
  \geq&t^2q^{(8d+1)/3}\left(t^2-2t\cdot5^{(-2\cdot2+2)/3}-5^{(-2-2)/3}-4t^2\cdot5^{-1}    -4t\cdot5^{(-2-2)/6}-5.03023t\right)\\
  \geq&t^2q^{(8d+1)/3}\left(0.2t^2-7.1t-0.2\right)\\
  \geq&t^2q^{(8d+1)/3}\left(\left(0.2\cdot36-7.1\right)\cdot36-0.2\right)>0.
  \end{split}
\end{equation}

\textbf{Case 2}: $d=3$ and $r\in\mathbb{F}_q^+$.

As in Subsection \ref{24zhijian}, it suffices to prove $|P_1(r)|>5.03023q^{(2d-2)/3}|\mathcal{E}|^3$.

Suppose $|\mathcal{E}|=tq^{(2d+1)/3}$ for some $t\geq9$.
Applying Theorem \ref{s1xiajie}, we have
  \begin{equation}
\begin{split}
  &|P_1(r)|-5.03023q^{(2d-2)/3}|\mathcal{E}|^3\\
  \geq&q^{-1}|\mathcal{E}|^{4}-2|\mathcal{E}|^{3}-2q^{d-1}|\mathcal{E}|^2-q^{-2}|\mathcal{E}|^4-2q^{(d-3)/2}|\mathcal{E}|^{3}-5.03023q^{(2d-2)/3}|\mathcal{E}|^3\\
  \geq&t^{4}q^{(8d+1)/3}      -2t^{3}q^{2d+1}    -2t^2q^{(7d-1)/3}       -t^4q^{(8d-2)/3}      -2t^{3}q^{(5d-1)/2}           -5.03023t^3q^{(8d+1)/3}\\
  =&t^2q^{(8d+1)/3}\left(t^2  -2tq^{(-2d+2)/3}   -2q^{(-d-2)/3}          -t^2q^{-1}            -2tq^{(-d-5)/6}               -5.03023t\right)\\
  \geq&t^2q^{(8d+1)/3}\left(t^2-2t\cdot3^{(-2\cdot3+2)/3}-2\cdot3^{(-3-2)/3}-t^2\cdot3^{-1}    -2t\cdot3^{(-3-5)/6}          -5.03023t\right)\\
  \geq&t^2q^{(8d+1)/3}\left(\frac{2}{3}t^2-5.95479t-0.3205\right)\\
  \geq&t^2q^{(8d+1)/3}\left(\frac{2}{3}\cdot81-5.95479\cdot9-0.3205\right)>0.
  \end{split}
\end{equation}

\textbf{Case 3}: $d=5$ and $r\in\mathbb{F}_q^+$.

As in Subsection \ref{ddengyu5}, it suffices to prove that $|\mathcal{E}|\geq12q^3$ implies $|P_1(r)|>1.34609|\mathcal{E}|^{11/3}$.

Suppose $|\mathcal{E}|=tq^{3}$ for some $t\geq12$.
Applying Theorem \ref{s1xiajie}, we have
  \begin{equation}
\begin{split}
  &|P_1(r)|-1.34609|\mathcal{E}|^{11/3}\\
  \geq&q^{-1}|\mathcal{E}|^{4}-2|\mathcal{E}|^{3}-2q^{d-1}|\mathcal{E}|^2-q^{-2}|\mathcal{E}|^4-2q^{(d-3)/2}|\mathcal{E}|^{3}-1.34609|\mathcal{E}|^{11/3}\\
  \geq&t^{4}q^{11}            -2t^{3}q^{9}       -2t^2q^{10}             -t^4q^{10}            -2t^{3}q^{10}                 -1.34609t^{11/3}q^{11}\\
  =&t^2q^{11}\left(t^2        -2tq^{-2}          -2q^{-1}                -t^2q^{-1}            -2tq^{-1}                     -1.34609t^{5/3}\right)\\
  \geq&t^2q^{11}\left(t^2        -\frac29t          -\frac23                -\frac{t^2}{3}    -\frac23t                     -1.34609t^{5/3}\right)\\
  =&t^2q^{11}\left(\frac{2}{3}t^2        -\frac89t          -\frac23                                   -1.34609t^{5/3}\right)\\
  \geq&t^2q^{11}\left(\frac{2}{3}t^2                -\frac{17}{18}t                                   -1.34609t^{5/3}\right)\\
    \geq&t^3q^{11}\left(\frac{2}{3}t               -\frac{17}{18}                            -1.34609t^{2/3}\right)\\
      =&t^3q^{11}\left(t^{2/3}\left(\frac{2}{3}t^{1/3}-1.34609\right)                -\frac{17}{18}\right)\\
       \geq&t^3q^{11}\left(12^{2/3}\left(\frac{2}{3}\cdot12^{1/3}-1.34609\right)                -\frac{17}{18}\right)>0.
  \end{split}
\end{equation}

\textbf{Case 4}: $q\geq5$, $d\geq6$ is even and $r\in\mathbb{F}_q^*$.

As in Subsection \ref{ddayudengyu6}, it suffices to prove that $|\mathcal{E}|\geq313q^{d/2}$ implies $|P_1(r)|>1.26844|\mathcal{E}|^{11/3}$.

Suppose $|\mathcal{E}|=tq^{d/2}$ for some $t\geq313$.
Applying Theorem \ref{s1xiajie}, we have
  \begin{equation}
\begin{split}
  &|P_1(r)|-1.26844|\mathcal{E}|^{11/3}\\
\geq&q^{-1}|\mathcal{E}|^{4}-2|\mathcal{E}|^{3}-q^{d-1}|\mathcal{E}|^2-4q^{-2}|\mathcal{E}|^4-4q^{(d-2)/2}|\mathcal{E}|^{3}-1.26844|\mathcal{E}|^{11/3}\\
\geq&t^{4}q^{2d-1}          -2t^{3}q^{3d/2}    -t^2q^{2d-1}           -4t^4q^{2d-2}          -4t^{3}q^{2d-1}               -1.26844t^{11/3}q^{11d/6}\\
  =&t^2q^{2d-1}\left(t^2    -2tq^{(2-d)/2}     -1                     -4t^2q^{-1}             -4t                          -1.26844t^{5/3}q^{(6-d)/6}\right)\\
\geq&t^2q^{2d-1}\left(t^2    -2t\cdot5^{(2-6)/2}     -1                     -4t^2\cdot5^{-1}             -4t                          -1.26844t^{5/3}\cdot5^{(6-6)/6}\right)\\
=&t^2q^{2d-1}\left(\frac{1}{5}t^2    -\frac{102}{25}t             -1                                                    -1.26844t^{5/3}\right)\\
\geq&t^2q^{2d-1}\left(\frac{1}{5}t^2    -\frac{49}{12}t                                                           -1.26844t^{5/3}\right)\\
=&t^3q^{2d-1}\left(t^{2/3}\left(\frac{1}{5}t^{1/3}    -1.26844\right)-\frac{49}{12}  \right)\\
\geq&t^3q^{2d-1}\left(313^{2/3}\left(\frac{1}{5}\cdot313^{1/3}    -1.26844\right)-\frac{49}{12}  \right)>0.
  \end{split}
\end{equation}

\textbf{Case 5}: $d\geq7$ is odd and $r\in\mathbb{F}_q^+$.

As in Subsection \ref{ddayudengyu6}, it suffices to prove that $|\mathcal{E}|\geq313q^{d/2}$ implies $|P_1(r)|>1.26844|\mathcal{E}|^{11/3}$.

Suppose $|\mathcal{E}|=tq^{d/2}$ for some $t\geq313$.
Applying Theorem \ref{s1xiajie}, we have
  \begin{equation}
\begin{split}
  &|P_1(r)|-1.26844|\mathcal{E}|^{11/3}\\
\geq&q^{-1}|\mathcal{E}|^{4}-2|\mathcal{E}|^{3}-2q^{d-1}|\mathcal{E}|^2-q^{-2}|\mathcal{E}|^4-2q^{(d-3)/2}|\mathcal{E}|^{3}-1.26844|\mathcal{E}|^{11/3}\\
\geq&t^{4}q^{2d-1}          -2t^{3}q^{3d/2}    -2t^2q^{2d-1}             -t^4q^{2d-2}            -2t^{3}q^{(4d-3)/2}                 -1.26844t^{11/3}q^{11d/6}\\
  =&t^2q^{2d-1}\left(t^2        -2tq^{(2-d)/2}          -2                -t^2q^{-1}            -2tq^{-1/2}                     -1.26844t^{5/3}q^{(6-d)/6}\right)\\
\geq&t^2q^{2d-1}\left(t^2        -2t\cdot3^{(2-7)/2}          -2                -t^2\cdot3^{-1}            -2t\cdot3^{-1/2}     -1.26844t^{5/3}\cdot3^{(6-7)/6}\right)\\
\geq&t^2q^{2d-1}\left(\frac{2}{3}t^2     -1.2894t                                  -1.05621t^{5/3}\right)\\
=&t^3q^{2d-1}\left(t^{2/3}\left(\frac{2}{3}t^{1/3}     -1.05621\right)-1.2894        \right)\\
\geq&t^3q^{2d-1}\left(313^{2/3}\left(\frac{2}{3}\cdot313^{1/3}     -1.05621\right)-1.2894        \right)>0.
  \end{split}
\end{equation}
\end{proof}

\section{Further remarks}
In \cite{MR4204716}, Greenleaf et al. considered the similar configurations problems in $\mathbb{R}^d$ and proved that if a compact set $\mathcal{E}$ has Hausdorff dimension greater than a threshold $s_{k, d}$, then there exist many pairs of $(k+1)$-point configurations which are similar by the scaling factor $r$. Rakhmonov \cite{2022arXiv220811579R} considered the similar $(d+1)$-simplices in $\mathbb{F}_p^d$ as well.

Recall that a $2$-star is a $2$-path. So it remains open for the problem about the similar $3$-paths in $\mathcal{E}\subseteq\mathbb{F}_q^d$. Unfortunately, it seems that there do not exist an appropriate inequality for $w_3(G)$ and $w_2(G)$.

\bibliographystyle{abbrv}
\bibliography{similarconfig_REF}

\begin{thebibliography}{10}

\bibitem{MR3592595}
M.~Bennett, D.~Hart, A.~Iosevich, J.~Pakianathan, and M.~Rudnev.
\newblock Group actions and geometric combinatorics in {$\Bbb{F}_q^d$}.
\newblock {\em Forum Math.}, 29(1):91--110, 2017.

\bibitem{MR184950}
G.~R. Blakley and P.~Roy.
\newblock A {H}\"{o}lder type inequality for symmetric matrices with
  nonnegative entries.
\newblock {\em Proc. Amer. Math. Soc.}, 16:1244--1245, 1965.

\bibitem{MR2917133}
J.~Chapman, M.~B. Erdo\u{g}an, D.~Hart, A.~Iosevich, and D.~Koh.
\newblock Pinned distance sets, {$k$}-simplices, {W}olff's exponent in finite
  fields and sum-product estimates.
\newblock {\em Math. Z.}, 271(1-2):63--93, 2012.

\bibitem{MR4201782}
X.~Du, L.~Guth, Y.~Ou, H.~Wang, B.~Wilson, and R.~Zhang.
\newblock Weighted restriction estimates and application to {F}alconer distance
  set problem.
\newblock {\em Amer. J. Math.}, 143(1):175--211, 2021.

\bibitem{MR4297185}
X.~Du, A.~Iosevich, Y.~Ou, H.~Wang, and R.~Zhang.
\newblock An improved result for {F}alconer's distance set problem in even
  dimensions.
\newblock {\em Math. Ann.}, 380(3-4):1215--1231, 2021.

\bibitem{MR4204716}
A.~Greenleaf, A.~Iosevich, and S.~Mkrtchyan.
\newblock Existence of similar point configurations in thin subsets of {$\Bbb
  R^d$}.
\newblock {\em Math. Z.}, 297(1-2):855--865, 2021.

\bibitem{MR4055179}
L.~Guth, A.~Iosevich, Y.~Ou, and H.~Wang.
\newblock On {F}alconer's distance set problem in the plane.
\newblock {\em Invent. Math.}, 219(3):779--830, 2020.

\bibitem{MR2775806}
D.~Hart, A.~Iosevich, D.~Koh, and M.~Rudnev.
\newblock Averages over hyperplanes, sum-product theory in vector spaces over
  finite fields and the {E}rd{\H{o}}s-{F}alconer distance conjecture.
\newblock {\em Trans. Amer. Math. Soc.}, 363(6):3255--3275, 2011.

\bibitem{MR3959878}
A.~Iosevich, D.~Koh, and H.~Parshall.
\newblock On the quotient set of the distance set.
\newblock {\em Mosc. J. Comb. Number Theory}, 8(2):103--115, 2019.

\bibitem{MR2336319}
A.~Iosevich and M.~Rudnev.
\newblock Erd{\H{o}}s distance problem in vector spaces over finite fields.
\newblock {\em Trans. Amer. Math. Soc.}, 359(12):6127--6142, 2007.

\bibitem{lagarias1983inequality}
J.~Lagarias, J.~Mazo, L.~Shepp, and B.~McKay.
\newblock An inequality for walks in a graph.
\newblock {\em SIAM Review}, 25(3):403--403, 1983.

\bibitem{MR1429394}
R.~Lidl and H.~Niederreiter.
\newblock {\em Finite fields}, volume~20 of {\em Encyclopedia of Mathematics
  and its Applications}.
\newblock Cambridge University Press, Cambridge, second edition, 1997.
\newblock With a foreword by P. M. Cohn.

\bibitem{2022arXiv220811579R}
F.~{Rakhmonov}.
\newblock {Distribution of similar configurations in subsets of
  $\mathbb{F}_q^d$}.
\newblock {\em arXiv e-prints}, page arXiv:2208.11579, Aug. 2022.

\end{thebibliography}
\end{document}